\documentclass{amsart}

\usepackage{amsmath,amsfonts,amssymb,amsmath,latexsym,multirow,color}

\newtheorem{theorem}{Theorem}[section]

\theoremstyle{definition}

\theoremstyle{remark}

\numberwithin{equation}{section}

\begin{document}

\title[Identities of symmetry for generalized Euler polynomials
]{$\begin{array}{c}
         \text{Identities of symmetry}\\
           \text{for generalized Euler polynomials}
       \end{array}$
 }

\author{dae san kim}
\address{Department of Mathematics, Sogang University, Seoul 121-742, Korea}
\curraddr{Department of Mathematics, Sogang University, Seoul
121-742, Korea} \email{dskim@sogong.ac.kr}

\begin{abstract}
In this paper, we derive eight basic identities of symmetry in three
variables related to generalized Euler polynomials and alternating
generalized power sums. All of these are new, since there have been
results only about identities of symmetry in two variables. The
derivations of identities are based on the $p$-adic fermionic
integral expression of the generating function for the generalized
Euler polynomials and the quotient of integrals that can be
expressed as the exponential generating function for the alternating
generalized power sums.

\end{abstract}


\date{}

\dedicatory{ }


\maketitle

\section{Introduction and preliminaries}
Let $p$ be a fixed odd prime. Throughout this paper, $\mathbb{Z}_p$,
$\mathbb{Q}_p$, $\mathbb{C}_p$ will respectively denote the ring of
$p$-adic integers, the field of $p$-adic rational numbers and the
completion of the algebraic closure of $\mathbb{Q}_p$. Let $d$ be a
fixed odd positive integer. Then we let
\begin{equation*}
X=X_d=\lim_{\overleftarrow{N}}\mathbb{Z}/dp^N \mathbb{Z},
\end{equation*}
and let $\pi:X\rightarrow \mathbb{Z}_p$ be the map given by the
inverse limit of the natural maps

\begin{equation*}
\mathbb{Z}/dp^N \mathbb{Z}\rightarrow \mathbb{Z}/p^N \mathbb{Z}.
\end{equation*}

\noindent If $g$ is a function on $\mathbb{Z}_p$, we will use the
same notation to denote the function $g\circ\pi$. Let
$\chi:(\mathbb{Z}/d\mathbb{Z})^*\rightarrow \overline{\mathbb{Q}}^*$
be a (primitive) Dirichlet character of conductor $d$. Then it will
be pulled back to $X$ via the natural map $X\rightarrow
\mathbb{Z}/d\mathbb{Z}$. Here we fix, once and for all, an imbedding
$\overline{\mathbb{Q}}\rightarrow \mathbb{C}_p$, so that $\chi$ is
regarded as a map of $X$ to $\mathbb{C}_p$ (cf. \cite{B8}).

For a continuous function $f:X\rightarrow \mathbb{C}_p$, the
$p$-adic fermionic integral of $f$ is defined by

\begin{equation*}
\int_{X}f(z)d\mu_{-1}(z)=\lim_{N\rightarrow\infty}\sum_{j=0}^{dp^N-1}f(j)(-1)^j.
\end{equation*}

\noindent Then it is easy to see that

\begin{equation}\label{N1}
\int_{X}f(z+1)d\mu_{-1}(z)+\int_{X}f(z)d\mu_{-1}(z)=2f(0).
\end{equation}

\noindent More generally, we deduce from (\ref{N1}) that, for any
odd positive integer $n$,

\begin{equation}\label{N2}
\int_{X}f(z+n)d\mu_{-1}(z)+\int_{X}f(z)d\mu_{-1}(z)=2\sum_{a=0}^{n-1}
(-1)^af(a),
\end{equation}

\noindent and that, for any even positive integer $n$,
\begin{equation*}
\int_{X}f(z+n)d\mu_{-1}(z)-\int_{X}f(z)d\mu_{-1}(z)=2\sum_{a=0}^{n-1}
(-1)^{a-1}f(a).
\end{equation*}

Let $|~|_p$ be the normalized absolute value of $\mathbb{C}_p$, such
that $|p|_p=\frac{1}{p}$, and let
\begin{equation}\label{N3}
E=\left\{ t\in \mathbb{C}_p ~\Big|~
|t|_p<p^{-\frac{1}{p-1}}\right\}.
\end{equation}

\noindent Then, for each fixed $t\in E$, the function $e^{zt}$ is
analytic on $\mathbb{Z}_p$ and hence considered as a function on
$X$, and, by applying (\ref{N2}) to $f$ with $f(z)=\chi(z)e^{zt}$,
we get the $p$-adic integral expression of the generating function
for the generalized Euler numbers $E_{n,\chi}$ attached to $\chi$:

\begin{equation}\label{N4}
\int_{X}\chi(z)e^{zt}d\mu_{-1}(z)=\frac{2}{e^{dt}+1}\sum_{a=0}^{d-1}(-1)^a
\chi(a)e^{at}=\sum_{n=0}^{\infty}E_{n,\chi}\frac{t^n}{n!} ~(t\in E).
\end{equation}

\noindent So we have the following $p$-adic integral expression of
the generating function for the generalized Euler polynomials
$E_{n,\chi}(x)$ attached to $\chi$:

\begin{equation}\label{N5}
\int_{X}\chi(z)e^{(x+z)t}d\mu_{-1}(z)=\frac{2e^{xt}}{e^{dt}+1}\sum_{a=0}^{d-1}(-1)^a
\chi(a)e^{at}=\sum_{n=0}^{\infty}E_{n,\chi}(x)\frac{t^n}{n!} ~(t\in
E,~x\in \mathbb{Z}_p).
\end{equation}

\noindent Also, from (\ref{N1}) we have:

\begin{equation}\label{N6}
\int_{X}e^{zt}d\mu_{-1}(z)=\frac{2}{e^{t}+1}~(t\in E).
\end{equation}

Let $T_k(n,\chi)$ denote the $k$th alternating generalized power sum
of the first $n+1$ nonnegative integers attached to $\chi$,  namely

\begin{equation}\label{N7}
T_k(n,\chi)=\sum_{a=0}^{n}(-1)^a\chi(a)a^k
=(-1)^0\chi(0)0^k+(-1)^1\chi(1)1^k+\dots+(-1)^n\chi(n)n^k.
\end{equation}

\noindent From (\ref{N4}), (\ref{N6}), and (\ref{N7}), one easily
derives the following identities: for any odd positive integer $w$,

\begin{align}
\frac
  {\int_{X}\chi(x)e^{xt}d\mu_{-1}(x)}{\int_{X}e^{wdyt}d\mu_{-1}(y)}
&=\frac{e^{wdt}+1}{e^{dt}+1}\sum_{a=0}^{d-1}(-1)^a\chi(a)e^{at}\\
&=\sum_{a=0}^{wd-1}(-1)^a \chi(a)e^{at}\\
&=\sum_{k=0}^{\infty}T_k(wd-1,\chi)\frac{t^k}{k!}~(t\in E).
\end{align}

\noindent In what follows, we will always assume that the $p$-adic
integrals of the various (twisted) exponential functions on $X$ are
defined for $t\in E$ (cf. (\ref{N3})), and therefore it will not be
mentioned.

\cite{B1}, \cite{B2}, \cite{B6}, \cite{B9} and \cite{B10} are some
of the previous works on identities of symmetry in two variables
involving Bernoulli polynomials and power sums. On the other hand,
for the first time we were able to produce in \cite{B4}  some
identities of symmetry in three variables related to Bernoulli
polynomials and power sums and to extend in \cite{B3} to the case of
generalized Bernoulli polynomials and generalized power sums.  Also,
\cite{B6} is about identities of symmetry in two variables for Euler
polynomials and alternating power sums and \cite{B5} is about those
in three variables for them.

In this paper, we will be able to produce 8 identities of symmetry
in three variables regarding generalized Euler polynomials and
alternating generalized power sums. The case of two variables was
treated in \cite{B7}.

The following is stated as Theorem \ref{T2} and  an example of the
full six symmetries in $w_1$, $w_2$, $w_3$.

\begin{align*}
&\sum_{k+l+m=n}{\binom{n}{k,l,m}}E_{k,\chi}(w_1y_1)E_{l,\chi}(w_2y_2)
T_{m}(w_3d-1,\chi)w_{1}^{l+m}w_{2}^{k+m}w_{3}^{k+l}\\
=&\sum_{k+l+m=n}{\binom{n}{k,l,m}}E_{k,\chi}(w_1y_1)E_{l,\chi}(w_3y_2)
T_{m}(w_2d-1,\chi)w_{1}^{l+m}w_{3}^{k+m}w_{2}^{k+l}\\
=&\sum_{k+l+m=n}{\binom{n}{k,l,m}}E_{k,\chi}(w_2y_1)E_{l,\chi}(w_1y_2)
T_{m}(w_3d-1,\chi)w_{2}^{l+m}w_{1}^{k+m}w_{3}^{k+l}\\
=&\sum_{k+l+m=n}{\binom{n}{k,l,m}}E_{k,\chi}(w_2y_1)E_{l,\chi}(w_3y_2)
T_{m}(w_1d-1,\chi)w_{2}^{l+m}w_{3}^{k+m}w_{1}^{k+l}\\
=&\sum_{k+l+m=n}{\binom{n}{k,l,m}}E_{k,\chi}(w_3y_1)E_{l,\chi}(w_2y_2)
T_{m}(w_1d-1,\chi)w_{3}^{l+m}w_{2}^{k+m}w_{1}^{k+l}\\
=&\sum_{k+l+m=n}{\binom{n}{k,l,m}}E_{k,\chi}(w_3y_1)E_{l,\chi}(w_1y_2)
T_{m}(w_2d-1,\chi)w_{3}^{l+m}w_{1}^{k+m}w_{2}^{k+l}.
\end{align*}

The derivations of identities are based on the $p$-adic integral
expression of the generating function for the generalized Euler
polynomials in (\ref{N5}) and the quotient of integrals in
(1.8)-(1.10) that can be expressed as the exponential generating
function for the alternating generalized power sums. These abundance
of symmetries would not be unearthed if such $p$-adic integral
representations had not been available. We indebted this idea to the
paper \cite{B7}.

\section{Several types of quotients of $p$-adic fermionic integrals}
Here we will introduce several types of quotients of $p$-adic
fermionic integrals on $X$ or $X^3$ from which some interesting
identities follow owing to the built-in symmetries in $w_1$, $w_2$,
$w_3$. In the following, $w_1$, $w_2$, $w_3$ are all positive
integers and all of the explicit expressions of integrals in
(\ref{N12}), (\ref{N14}), (\ref{N16}), and (\ref{N18}) are obtained
from the identities in (\ref{N4}) and (\ref{N6}). To ease notations,
from now on we will suppress $\mu_{-1}$ and denote, for example,
$d\mu_{-1}(x)$ simply by $dx$.

\noindent (a) Type $\Lambda_{23}^{i}$ (for $i=0,1,2,3$)

\begin{equation}\label{N11}
I(\Lambda_{23}^{i})=\frac{\int_{X^3}\chi(x_1)\chi(x_2)\chi(x_3)
e^{(w_2w_3x_1+w_1w_3x_2+w_1w_2x_3+w_1w_2w_3(\sum_{j=1}^{3-i}y_j))t}dx_1dx_2dx_3}
 {(\int_{X}e^{dw_1w_2w_3x_4t}dx_4)^i}\\
\end{equation}
\begin{equation}\label{N12}
\begin{split}
=&\frac{2^{3-i}e^{w_1w_2w_3(\sum_{j=1}^{3-i}y_j)t}(e^{dw_1w_2w_3t}+1)^i}
{(e^{dw_2w_3t}+1)(e^{dw_1w_3t}+1)(e^{dw_1w_2t}+1)}\\
&\times(\sum_{a=0}^{d-1}(-1)^a\chi(a)e^{aw_2w_3t})(\sum_{a=0}^{d-1}(-1)^a\chi(a)e^{aw_1w_3t})
(\sum_{a=0}^{d-1}(-1)^a\chi(a)e^{aw_1w_2t}).
\end{split}
\end{equation}

\noindent (b) Type $\Lambda_{13}^{i}$ (for $i=0,1,2,3$)

\begin{equation}\label{N13}
I(\Lambda_{13}^{i})=\frac{\int_{X^3}\chi(x_1)\chi(x_2)\chi(x_3)
e^{(w_1x_1+w_2x_2+w_3x_3+w_1w_2w_3(\sum_{j=1}^{3-i}y_j))t}dx_1dx_2dx_3}
 {(\int_{X}e^{dw_1w_2w_3x_4t}dx_4)^i}\\
\end{equation}
\begin{equation}\label{N14}
\begin{split}
=&\frac{2^{3-i}e^{w_1w_2w_3(\sum_{j=1}^{3-i}y_j)t}(e^{dw_1w_2w_3t}+1)^i}
{(e^{dw_1t}+1)(e^{dw_2t}+1)(e^{dw_3t}+1)}\\
&\times(\sum_{a=0}^{d-1}(-1)^a\chi(a)e^{aw_1t})(\sum_{a=0}^{d-1}(-1)^a\chi(a)e^{aw_2t})
(\sum_{a=0}^{d-1}(-1)^a\chi(a)e^{aw_3t}).
\end{split}
\end{equation}

\noindent (c-0) Type $\Lambda_{12}^{0}$

\begin{equation}\label{N15}
I(\Lambda_{12}^{0})=\int_{X^3}\chi(x_1)\chi(x_2)\chi(x_3)
e^{(w_1x_1+w_2x_2+w_3w_3+w_2w_3y+w_1w_3y+w_1w_2y)t}dx_1dx_2dx_3
\end{equation}
\begin{equation}\label{N16}
\begin{split}
=&\frac{8e^{(w_2w_3+w_1w_3+w_1w_2)yt}}
{(e^{dw_1t}+1)(e^{dw_2t}+1)(e^{dw_3t}+1)}\\
&\times(\sum_{a=0}^{d-1}(-1)^a\chi(a)e^{aw_1t})(\sum_{a=0}^{d-1}(-1)^a\chi(a)e^{aw_2t})
(\sum_{a=0}^{d-1}(-1)^a\chi(a)e^{aw_3t}).
\end{split}
\end{equation}

\noindent (c-1) Type $\Lambda_{12}^{1}$

\begin{equation}\label{N17}
I(\Lambda_{12}^{1})=\frac{\int_{X^3}\chi(x_1)\chi(x_2)\chi(x_3)
e^{(w_1x_1+w_2x_2+w_3x_3)t}dx_1dx_2dx_3}
 {\int_{X^3}e^{d(w_2w_3z_1+w_1w_3z_2+w_1w_2z_3)t}dz_1z_2z_3}\\
\end{equation}
\begin{equation}\label{N18}
\begin{split}
=&\frac{(e^{dw_2w_3t}+1)(e^{dw_1w_3t}+1)(e^{dw_1w_2t}+1)}
{(e^{dw_1t}+1)(e^{dw_2t}+1)(e^{dw_3t}+1)}\\
&\times(\sum_{a=0}^{d-1}(-1)^a\chi(a)e^{aw_1t})
(\sum_{a=0}^{d-1}(-1)^a\chi(a)e^{aw_2t})
(\sum_{a=0}^{d-1}(-1)^a\chi(a)e^{aw_3t}).
\end{split}
\end{equation}

All of the above $p$-adic integrals of various types are invariant
under all permutations of $w_1$, $w_2$, $w_3$, as one can see either
from $p$-adic integral representations in (\ref{N11}), (\ref{N13}),
(\ref{N15}), and (\ref{N17}) or from their explicit evaluations in
(\ref{N12}), (\ref{N14}), (\ref{N16}), and (\ref{N18}).

\section{Identities for generalized Euler polynomials}

In the following, $w_1$, $w_2$, $w_3$ are all odd positive integers
except for (a-0) and (c-0), where they are any positive integers.
First, let's consider Type $\Lambda_{23}^{i}$, for each $i=0,1,2,3$.
The following results can be easily obtained from (\ref{N5}) and
(1.8)-(1.10).

\noindent (a-0)
\begin{equation}\label{N19}
\begin{split}
I(&\Lambda_{23}^{0})\\
&=\int_{X}\chi(x_1)e^{w_2w_3(x_1+w_1y_1)t}dx_1
\int_{X}\chi(x_2)e^{w_1w_3(x_2+w_2y_2)t}dx_2\\
&\qquad\qquad\qquad\qquad\qquad\qquad\qquad\qquad\times
\int_{X}\chi(x_3)e^{w_1w_2(x_3+w_3y_3)t}dx_3\\
&=(\sum_{k=0}^{\infty}\frac{E_{k,\chi}(w_1y_1)}{k!}(w_2w_3t)^k)
(\sum_{l=0}^{\infty}\frac{E_{l,\chi}(w_2y_2)}{l!}(w_1w_3t)^l)\\
&\qquad\qquad\qquad\qquad\qquad\qquad\qquad\qquad\times
(\sum_{m=0}^{\infty}\frac{E_{m,\chi}(w_3y_3)}{m!}(w_1w_2t)^m)\\
&=\sum_{n=0}^{\infty}(\sum_{k+l+m=n}\binom{n}{k,l,m}E_{k,\chi}(w_1y_1)E_{l,\chi}(w_2y_2)
E_{m,\chi}(w_3y_3)\\
&\qquad\qquad\qquad\qquad\qquad\qquad\qquad\qquad\qquad\times
w_{1}^{l+m}w_{2}^{k+m}w_{3}^{k+l})\frac{t^n}{n!},
\end{split}
\end{equation}

\noindent where the inner sum is over all nonnegative integers
$k,~l,~m$ with $k+l+m=n$,and

\begin{equation}\label{N20}
\binom{n}{k,l,m}=\frac{n!}{k!l!m!}.
\end{equation}

\noindent (a-1) Here we write $I(\Lambda_{23}^{1})$ in two different
ways:

\begin{equation}\label{N21}
\begin{split}
\noindent
(1)~I(\Lambda_{23}^{1})&=\int_{X}\chi(x_1)e^{w_2w_3(x_1+w_1y_1)t}dx_1
\int_{X}\chi(x_2)e^{w_1w_3(x_2+w_2y_2)t}dx_2\\
&\qquad\qquad\qquad\qquad\qquad\qquad\qquad
\times\frac{\int_{X}\chi(x_3)e^{w_1w_2x_3t}dx_3}{\int_{X}e^{dw_1w_2w_3x_4t}dx_4}\\
&=(\sum_{k=0}^{\infty}E_{k,\chi}(w_1y_1)\frac{(w_2w_3t)^k}{k!})
(\sum_{l=0}^{\infty}E_{l,\chi}(w_2y_2)\frac{(w_1w_3t)^l}{l!})\\
&\qquad\qquad\qquad\qquad\qquad\qquad\qquad\quad
\times(T_m(w_3d-1,\chi)\frac{(w_1w_2t)^m}{m!})\\
\end{split}
\end{equation}

\begin{equation}\label{N22}
\begin{split}
=&\sum_{n=0}^{\infty}(\sum_{k+l+m=n}\binom{n}{k,l,m}E_{k,\chi}(w_1y_1)E_{l,\chi}(w_2y_2)\\
&\qquad\qquad\qquad\qquad\quad\times
T_m(w_3d-1,\chi)w_{1}^{l+m}w_{2}^{k+m}w_{3}^{k+l})\frac{t^n}{n!}.
\end{split}
\end{equation}

\noindent (2) Invoking (1.9), (\ref{N21}) can also be written as
\begin{equation}\label{N23}
\begin{split}
I(&\Lambda_{23}^{1})\\
&=\sum_{a=0}^{w_3d-1}(-1)^a\chi(a)
\int_{X}\chi(x_1)e^{w_2w_3(x_1+w_1y_1)t}dx_1\\
&\qquad\qquad\qquad\qquad\qquad\qquad\qquad\qquad\times
\int_{X}\chi(x_2)e^{w_1w_3(x_2+w_2y_2+\frac{w_2}{w_3}a)t}dx_2\\
&=\sum_{a=0}^{w_3d-1}(-1)^a\chi(a)
(\sum_{k=0}^{\infty}E_{k,\chi}(w_1y_1)\frac{(w_2y_3t)^k}{k!})\\
&\qquad\qquad\qquad\qquad\qquad\qquad\qquad\qquad\times
(\sum_{l=0}^{\infty}E_{l,\chi}(w_2y_2+\frac{w_2}{w_3}a)\frac{(w_1y_3t)^l}{l!})\\
&=\sum_{n=0}^{\infty}(w_3^n\sum_{k=0}^{n}\binom{n}{k}E_{k,\chi}(w_1y_1)
\sum_{a=0}^{w_3d-1}(-1)^a\chi(a)\\
&\qquad\qquad\qquad\qquad\qquad\qquad\qquad\qquad\times E_{n-k,\chi}
(w_2y_2+\frac{w_2}{w_3}a)w_1^{n-k}
w_2^k)\frac{t^n}{n!}.\\
&
\end{split}
\end{equation}

\noindent (a-2) Here we write $I(\Lambda_{23}^{2})$ in three
different ways:

\begin{equation}\label{N24}
\begin{split}
\noindent(1)~I(&\Lambda_{23}^{2})\\
&=\int_{X}\chi(x_1)e^{w_2w_3(x_1+w_1y_1)t}dx_1
\times\frac{\int_{X}\chi(x_2)e^{w_1w_3x_2t}dx_2}{\int_{X}e^{dw_1w_2w_3x_4t}dx_4}\\
&\qquad\qquad\qquad\qquad\qquad\qquad\qquad\qquad\qquad\times\frac{\int_{X}\chi(x_3)e^{w_1w_2x_3t}dx_3}{\int_{X}e^{dw_1w_2w_3x_4t}dx_4}\\
&=(\sum_{k=0}^{\infty}E_{k,\chi}(w_1y_1)\frac{(w_2w_3t)^k}{k!})
(\sum_{l=0}^{\infty}T_{l}(w_2d-1,\chi)\frac{(w_1w_3t)^l}{l!})\\
&\qquad\qquad\qquad\qquad\qquad\qquad\times
(\sum_{m=0}^{\infty}T_{m}(w_3d-1,\chi)\frac{(w_1w_2t)^m}{m!})
\\
\end{split}
\end{equation}

\begin{equation}\label{N25}
\begin{split}
=&\sum_{n=0}^{\infty}(\sum_{k+l+m=n}\binom{n}{k,l,m}E_{k,\chi}(w_1y_1)T_l(w_2d-1,\chi)\\
&\qquad\qquad\qquad\qquad\qquad\qquad\quad\times
T_m(w_3d-1,\chi)w_{1}^{l+m}w_{2}^{k+m}w_{3}^{k+l})\frac{t^n}{n!}.
\end{split}
\end{equation}

\noindent (2) Invoking (1.9), (\ref{N24}) can also be written as

\begin{equation}\label{N26}
\begin{split}
I(&\Lambda_{23}^{2})=\sum_{a=0}^{w_2d-1}(-1)^a\chi(a)
\int_{X}\chi(x_1)e^{w_2w_3(x_1+w_1y_1+\frac{w_1}{w_2}a)t}dx_1\\
&\qquad\qquad\qquad\qquad\qquad\qquad\qquad\qquad\qquad
\times\frac{\int_{X}\chi(x_3)e^{w_1w_2x_3t}dx_3}{\int_{X}e^{dw_1w_2w_3x_4t}dx_4}
\end{split}
\end{equation}
\begin{equation*}
\begin{split}
&\qquad\qquad=\sum_{a=0}^{w_2d-1}(-1)^a\chi(a)
(\sum_{k=0}^{\infty}E_{k,\chi}(w_1y_1+\frac{w_1}{w_2}a)\frac{(w_2w_3t)^k}{k!})\\
&\qquad\qquad\qquad\qquad\qquad\qquad\qquad\qquad\times(\sum_{l=0}^{\infty}T_{l}(w_3d-1,\chi)\frac{(w_1w_2t)^l}{l!})\\
\end{split}
\end{equation*}

\begin{equation}\label{N27}
\begin{split}
&=\sum_{n=0}^{\infty}(w_2^n\sum_{k=0}^{n}\binom{n}{k}
\sum_{a=0}^{w_2d-1}(-1)^a\chi(a)
E_{k,\chi}(w_1y_1+\frac{w_1}{w_2}a)\\
&\qquad\qquad\qquad\qquad\qquad\qquad\qquad\qquad\times
T_{n-k}(w_3d-1,\chi)w_{1}^{n-k}w_3^{k})\frac{t^n}{n!}.
\end{split}
\end{equation}

\noindent (3) Invoking (1.9) once again, (\ref{N26}) can be written
as

\begin{equation*}
\begin{split}
I(&\Lambda_{23}^{2})\\
&=\sum_{a=0}^{w_2d-1}(-1)^a\chi(a) \sum_{b=0}^{w_3d-1}(-1)^b\chi(b)
\int_{X}\chi(x_1)e^{w_2w_3(x_1+w_1y_1+\frac{w_1}{w_2}a+\frac{w_1}{w_3}b)t}dx_1\\
&=\sum_{a=0}^{w_2d-1}(-1)^a\chi(a) \sum_{b=0}^{w_3d-1}(-1)^b\chi(b)
\sum_{n=0}^{\infty}E_{n,\chi}(w_1y_1+\frac{w_1}{w_2}a+\frac{w_1}{w_3}b)
\frac{(w_2w_3t)^n}{n!}
\end{split}
\end{equation*}

\begin{equation}\label{N28}
=\sum_{n=0}^{\infty}((w_2w_3)^n\sum_{a=0}^{w_2d-1}\sum_{b=0}^{w_3d-1}(-1)^{a+b}\chi(ab)
E_{n,\chi}(w_1y_1+\frac{w_1}{w_2}a+\frac{w_1}{w_3}b))
\frac{t^n}{n!}.\qquad
\end{equation}

\noindent (a-3)
\begin{equation*}
\begin{split}
I(&\Lambda_{23}^{3})\\
&=\frac{\int_{X}\chi(x_1)e^{w_2w_3x_1t}dx_1}{\int_{X}e^{dw_1w_2w_3x_4t}dx_4}
\times\frac{\int_{X}\chi(x_2)e^{w_1w_3x_2t}dx_2}{\int_{X}e^{dw_1w_2w_3x_4t}dx_4}
\times\frac{\int_{X}\chi(x_3)e^{w_1w_2x_3t}dx_3}{\int_{X}e^{dw_1w_2w_3x_4t}dx_4}\\
&=(\sum_{k=0}^{\infty}T_{k}(w_1d-1,\chi)\frac{(w_2w_3t)^k}{k!})
(\sum_{l=0}^{\infty}T_{l}(w_2d-1,\chi)\frac{(w_1w_3t)^l}{l!})\\
&\qquad\qquad\qquad\qquad\qquad\qquad\qquad\qquad\times
(\sum_{m=0}^{\infty}T_{m}(w_3d-1,\chi)\frac{(w_1w_2t)^m}{m!})
\\
\end{split}
\end{equation*}

\begin{equation}\label{N29}
\begin{split}
=&\sum_{n=0}^{\infty}(\sum_{k+l+m=n}\binom{n}{k,l,m}T_k(w_1d-1,\chi)
T_l(w_2d-1,\chi)T_m(w_3d-1,\chi)\qquad\\
&\qquad\qquad\qquad\qquad\qquad\qquad\quad\times
w_{1}^{l+m}w_{2}^{k+m}w_{3}^{k+l})\frac{t^n}{n!}.\\
\end{split}
\end{equation}

\noindent (b) For Type $\Lambda_{13}^{i}~($i=0,1,2,3$)$, we may
consider the analogous things to the ones in (a-0), (a-1), (a-2),
and (a-3). However, these do not lead us to new identities. Indeed,
if we substitute $w_2w_3$, $w_1w_3$, $w_1w_2$ respectively for
$w_1$, $w_2$, $w_3$ in (\ref{N11}), this amounts to replacing $t$ by
$w_1w_2w_3t$ in (\ref{N13}). So, upon replacing $w_1$, $w_2$, $w_3$
respectively by $w_2w_3$, $w_1w_3$, $w_1w_2$ and dividing by $(w_1
w_2 w_3)^n$, in each of the expressions of (\ref{N19}), (\ref{N22}),
(\ref{N23}), (\ref{N25}), (\ref{N27})-(\ref{N29}), we will get the
corresponding symmetric identities for Type $\Lambda_{13}^{i}$
($i=0,1,2,3$).

\noindent (c-0)
\begin{equation*}
\begin{split}
I(&\Lambda_{12}^{0})\\
&=\int_{X}\chi(x_1)e^{w_1(x_1+w_2y)t}dx_1
\int_{X}\chi(x_2)e^{w_2(x_2+w_3y)t}dx_2\\
&\qquad\qquad\qquad\qquad\qquad\qquad\qquad
\times\int_{X}\chi(x_3)e^{w_3(x_3+w_1y)t}dx_3\\
&=(\sum_{k=0}^{\infty}\frac{E_{k,\chi}(w_2y)}{k!}(w_1t)^k)
(\sum_{l=0}^{\infty}\frac{E_{l,\chi}(w_3y)}{l!}(w_2t)^l)
(\sum_{m=0}^{\infty}\frac{E_{m,\chi}(w_1y)}{m!}(w_3t)^m)
\end{split}
\end{equation*}

\begin{equation}\label{N30}
=\sum_{n=0}^{\infty}(\sum_{k+l+m=n}\binom{n}{k,l,m} E_{k,\chi}(w_2y)
E_{l,\chi}(w_3y)E_{m,\chi}(w_1y)w_{1}^{k}w_{2}^{l}w_{3}^{m}
)\frac{t^n}{n!}.\qquad
\end{equation}

\noindent (c-1)
\begin{equation*}
\begin{split}
I(&\Lambda_{12}^{1})\\
&=\frac{\int_{X}\chi(x_1)e^{w_1x_1t}dx_1}{\int_{X}e^{dw_1w_2z_3t}dz_3}
\times\frac{\int_{X}\chi(x_2)e^{w_2x_2t}dx_2}{\int_{X}e^{dw_2w_3z_1t}dz_1}
\times\frac{\int_{X}\chi(x_3)e^{w_3x_3t}dx_3}{\int_{X}e^{dw_3w_1z_2t}dz_2}\\
&=(\sum_{k=0}^{\infty}T_{k}(w_2d-1,\chi)\frac{(w_1t)^k}{k!})
(\sum_{l=0}^{\infty}T_{l}(w_3d-1,\chi)\frac{(w_2t)^l}{l!})\\
&\qquad\qquad\qquad\qquad\qquad\qquad\qquad\qquad\times
(\sum_{m=0}^{\infty}T_{m}(w_1d-1,\chi)\frac{(w_3t)^m}{m!})
\\
\end{split}
\end{equation*}

\begin{equation}\label{N31}
\begin{split}
=&\sum_{n=0}^{\infty}(\sum_{k+l+m=n}\binom{n}{k,l,m}T_k(w_2d-1,\chi)
T_l(w_3d-1,\chi)\\
&\qquad\qquad\qquad\qquad\qquad\qquad\quad\times
T_m(w_1d-1,\chi)w_{1}^{k}w_{2}^{l}w_{3}^{m})\frac{t^n}{n!}.\qquad\quad
\end{split}
\end{equation}

\section{Main theorems}
As we noted earlier in the last paragraph of Section 2, the various
types of quotients of $p$-adic fermionic integrals are invariant
under any permutation of $w_1$, $w_2$, $w_3$. So the corresponding
expressions in Section 3 are also invariant under any permutation of
$w_1$, $w_2$, $w_3$. Thus our results about identities of symmetry
will be immediate consequences of this observation.

However, not all permutations of an expression in Section 3 yield
distinct ones. In fact, as these expressions are obtained by
permuting $w_1$, $w_2$, $w_3$ in a single one labeled by them, they
can be viewed as a group in a natural manner and hence it is
isomorphic to a quotient of $S_3$. In particular, the number of
possible distinct expressions are $1$, $2$, $3$, or $6$ (a-0),
(a-1(1)), (a-1(2)), and (a-2(2)) give the full six identities of
symmetry, (a-2(1)) and (a-2(3)) yield three identities of symmetry,
and (c-0) and (c-1) give two identities of symmetry, while the
expression in (a-3) yields no identities of symmetry.

Here we will just consider the cases of Theorems \ref{T4} and
\ref{T8}, leaving the others as easy exercises for the reader. As
for the case of Theorem \ref{T4}, in addition to (4.11)-(4.13), we
get the following three ones:
\begin{align}
&\sum_{k+l+m=n}\binom{n}{k,l,m}E_{k,\chi}(w_1y_1)T_l(w_3d-1,\chi)T_m(w_2d-1,\chi)
w_{1}^{l+m}w_{3}^{k+m}w_{2}^{k+l},\\
&=\sum_{k+l+m=n}\binom{n}{k,l,m}E_{k,\chi}(w_2y_1)T_l(w_1d-1,\chi)T_m(w_3d-1,\chi)
w_{2}^{l+m}w_{1}^{k+m}w_{3}^{k+l},\\
&=\sum_{k+l+m=n}\binom{n}{k,l,m}E_{k,\chi}(w_3y_1)T_l(w_2d-1,\chi)T_m(w_1d-1,\chi)
w_{3}^{l+m}w_{2}^{k+m}w_{1}^{k+l}.
\end{align}

\noindent But, by interchanging $l$ and $m$, we see that (4.1),
(4.2), and (4.3) are respectively equal to (4.11), (4.12), and
(4.13). As to Theorem \ref{T8}, in addition to (4.17) and (4.18), we
have:
\begin{align}
&\sum_{k+l+m=n}\binom{n}{k,l,m}T_k(w_2d-1,\chi)T_l(w_3d-1,\chi)T_m(w_1d-1,\chi)
w_{1}^{k}w_{2}^{l}w_{3}^{m},\\
&=\sum_{k+l+m=n}\binom{n}{k,l,m}T_k(w_3d-1,\chi)T_l(w_1d-1,\chi)T_m(w_2d-1,\chi)
w_{2}^{k}w_{3}^{l}w_{1}^{m},\\
&=\sum_{k+l+m=n}\binom{n}{k,l,m}T_k(w_3d-1,\chi)T_l(w_2d-1,\chi)T_m(w_1d-1,\chi)
w_{1}^{k}w_{3}^{l}w_{2}^{m},\\
&=\sum_{k+l+m=n}\binom{n}{k,l,m}T_k(w_2d-1,\chi)T_l(w_1d-1,\chi)T_m(w_3d-1,\chi)
w_{3}^{k}w_{2}^{l}w_{1}^{m}.
\end{align}

\noindent However, (4.4) and (4.5) are equal to (4.17), as we can
see by applying the permutations $k\rightarrow l$, $l\rightarrow m$,
$m\rightarrow k$ for (4.4) and  $k\rightarrow m$, $l \rightarrow k$,
$m\rightarrow l$ for (4.5). Similarly, we see that (4.6) and (4.7)
are equal to (4.18), by applying permutations $k\rightarrow l$,
$l\rightarrow m$, $m\rightarrow k$ for (4.6) and $k\rightarrow m$,
$l\rightarrow k$, $m\rightarrow l$ for (4.7).

\begin{theorem}\label{T1}
Let $w_1$, $w_2$, $w_3$ be any positive integers. Then we have:
\begin{equation}\label{N39}
\begin{split}
&\sum_{k+l+m=n}\binom{n}{k,l,m}E_{k,\chi}(w_1y_1)E_{l,\chi}(w_2y_2)E_{m,\chi}(w_3y_3)
w_{1}^{l+m}w_{2}^{k+m}w_{3}^{k+l}\\
=&\sum_{k+l+m=n}\binom{n}{k,l,m}E_{k,\chi}(w_1y_1)E_{l,\chi}(w_3y_2)E_{m,\chi}(w_2y_3)
w_{1}^{l+m}w_{3}^{k+m}w_{2}^{k+l}\\
=&\sum_{k+l+m=n}\binom{n}{k,l,m}E_{k,\chi}(w_2y_1)E_{l,\chi}(w_1y_2)E_{m,\chi}(w_3y_3)
w_{2}^{l+m}w_{1}^{k+m}w_{3}^{k+l}\\
=&\sum_{k+l+m=n}\binom{n}{k,l,m}E_{k,\chi}(w_2y_1)E_{l,\chi}(w_3y_2)E_{m,\chi}(w_1y_3)
w_{2}^{l+m}w_{3}^{k+m}w_{1}^{k+l}\\
=&\sum_{k+l+m=n}\binom{n}{k,l,m}E_{k,\chi}(w_3y_1)E_{l,\chi}(w_1y_2)E_{m,\chi}(w_2y_3)
w_{3}^{l+m}w_{1}^{k+m}w_{2}^{k+l}\\
=&\sum_{k+l+m=n}\binom{n}{k,l,m}E_{k,\chi}(w_3y_1)E_{l,\chi}(w_2y_2)E_{m,\chi}(w_1y_3)
w_{3}^{l+m}w_{2}^{k+m}w_{1}^{k+l}.
\end{split}
\end{equation}
\end{theorem}

\begin{theorem}\label{T2}
Let $w_1$, $w_2$, $w_3$ be any odd positive integers. Then we have:
\begin{equation}\label{N40}
\begin{split}
&\sum_{k+l+m=n}\binom{n}{k,l,m}E_{k,\chi}(w_1y_1)E_{l,\chi}(w_2y_2)T_m(w_3d-1,\chi)
w_{1}^{l+m}w_{2}^{k+m}w_{3}^{k+l}\\
=&\sum_{k+l+m=n}\binom{n}{k,l,m}E_{k,\chi}(w_1y_1)E_{l,\chi}(w_3y_2)T_m(w_2d-1,\chi)
w_{1}^{l+m}w_{3}^{k+m}w_{2}^{k+l}\\
=&\sum_{k+l+m=n}\binom{n}{k,l,m}E_{k,\chi}(w_2y_1)E_{l,\chi}(w_1y_2)T_m(w_3d-1,\chi)
w_{2}^{l+m}w_{1}^{k+m}w_{3}^{k+l}\\
=&\sum_{k+l+m=n}\binom{n}{k,l,m}E_{k,\chi}(w_2y_1)E_{l,\chi}(w_3y_2)T_m(w_1d-1,\chi)
w_{2}^{l+m}w_{3}^{k+m}w_{1}^{k+l}\\
=&\sum_{k+l+m=n}\binom{n}{k,l,m}E_{k,\chi}(w_3y_1)E_{l,\chi}(w_2y_2)T_m(w_1d-1,\chi)
w_{3}^{l+m}w_{2}^{k+m}w_{1}^{k+l}\\
=&\sum_{k+l+m=n}\binom{n}{k,l,m}E_{k,\chi}(w_3y_1)E_{l,\chi}(w_1y_2)T_m(w_2d-1,\chi)
w_{3}^{l+m}w_{1}^{k+m}w_{2}^{k+l}.
\end{split}
\end{equation}
\end{theorem}

\begin{theorem}\label{T3}
Let $w_1$, $w_2$, $w_3$ be any odd positive integers. Then we have:
\begin{equation}\label{N41}
\begin{split}
&w_1^n\sum_{k=0}^{n}\binom{n}{k}E_{k,\chi}(w_3y_1)\sum_{a=0}^{w_1d-1}(-1)^a
\chi(a)
E_{n-k,\chi}(w_2y_2+\frac{w_2}{w_1}a)w_3^{n-k}w_2^k\\
=&w_1^n\sum_{k=0}^{n}\binom{n}{k}E_{k,\chi}(w_2y_1)\sum_{a=0}^{w_1d-1}(-1)^a\chi(a)
E_{n-k,\chi}(w_3y_2+\frac{w_3}{w_1}a)w_2^{n-k}w_3^k\\
=&w_2^n\sum_{k=0}^{n}\binom{n}{k}E_{k,\chi}(w_3y_1)\sum_{a=0}^{w_2d-1}(-1)^a\chi(a)
E_{n-k,\chi}(w_1y_2+\frac{w_1}{w_2}a)w_3^{n-k}w_1^k\\
=&w_2^n\sum_{k=0}^{n}\binom{n}{k}E_{k,\chi}(w_1y_1)\sum_{a=0}^{w_2d-1}(-1)^a\chi(a)
E_{n-k,\chi}(w_3y_2+\frac{w_3}{w_2}a)w_1^{n-k}w_3^k
\end{split}
\end{equation}
\begin{equation*}
\begin{split}
=&w_3^n\sum_{k=0}^{n}\binom{n}{k}E_{k,\chi}(w_2y_1)\sum_{a=0}^{w_3d-1}(-1)^a\chi(a)
E_{n-k,\chi}(w_1y_2+\frac{w_1}{w_3}a)w_2^{n-k}w_1^k\\
=&w_3^n\sum_{k=0}^{n}\binom{n}{k}E_{k,\chi}(w_1y_1)\sum_{a=0}^{w_3d-1}(-1)^a\chi(a)
E_{n-k,\chi}(w_2y_2+\frac{w_2}{w_3}a)w_1^{n-k}w_2^k.
\end{split}
\end{equation*}
\end{theorem}

\begin{theorem}\label{T4}
Let $w_1$, $w_2$, $w_3$ be any odd positive integers. Then we have
the following three symmetries in $w_1$, $w_2$, $w_3$:
\begin{align}
&\sum_{k+l+m=n}\binom{n}{k,l,m}E_{k,\chi}(w_1y_1)T_l(w_2d-1,\chi)T_m(w_3d-1,\chi)
w_{1}^{l+m}w_{2}^{k+m}w_{3}^{k+l}\\
=&\sum_{k+l+m=n}\binom{n}{k,l,m}E_{k,\chi}(w_2y_1)T_l(w_3d-1,\chi)T_m(w_1d-1,\chi)
w_{2}^{l+m}w_{3}^{k+m}w_{1}^{k+l}\\
=&\sum_{k+l+m=n}\binom{n}{k,l,m}E_{k,\chi}(w_3y_1)T_l(w_1d-1,\chi)T_m(w_2d-1,\chi)
w_{3}^{l+m}w_{1}^{k+m}w_{2}^{k+l}.
\end{align}
\end{theorem}

\begin{theorem}\label{T5}
Let $w_1$, $w_2$, $w_3$ be any odd positive integers. Then we have:
\begin{equation}\label{N45}
\begin{split}
&w_1^n\sum_{k=0}^{n}\binom{n}{k}\sum_{a=0}^{w_1d-1}(-1)^a\chi(a)
E_{k,\chi}(w_2y_1+\frac{w_2}{w_1}a)T_{n-k}(w_3d-1,\chi)w_2^{n-k}w_3^k\\
=&w_1^n\sum_{k=0}^{n}\binom{n}{k}\sum_{a=0}^{w_1d-1}(-1)^a\chi(a)
E_{k,\chi}(w_3y_1+\frac{w_3}{w_1}a)T_{n-k}(w_2d-1,\chi)w_3^{n-k}w_2^k\\
=&w_2^n\sum_{k=0}^{n}\binom{n}{k}\sum_{a=0}^{w_2d-1}(-1)^a\chi(a)
E_{k,\chi}(w_1y_1+\frac{w_1}{w_2}a)T_{n-k}(w_3d-1,\chi)w_1^{n-k}w_3^k\\
=&w_2^n\sum_{k=0}^{n}\binom{n}{k}\sum_{a=0}^{w_2d-1}(-1)^a\chi(a)
E_{k,\chi}(w_3y_1+\frac{w_3}{w_2}a)T_{n-k}(w_1d-1,\chi)w_3^{n-k}w_1^k\\
=&w_3^n\sum_{k=0}^{n}\binom{n}{k}\sum_{a=0}^{w_3d-1}(-1)^a\chi(a)
E_{k,\chi}(w_1y_1+\frac{w_1}{w_3}a)T_{n-k}(w_2d-1,\chi)w_1^{n-k}w_2^k\\
=&w_3^n\sum_{k=0}^{n}\binom{n}{k}\sum_{a=0}^{w_3d-1}(-1)^a\chi(a)
E_{k,\chi}(w_2y_1+\frac{w_2}{w_3}a)T_{n-k}(w_1d-1,\chi)w_2^{n-k}w_1^k.
\end{split}
\end{equation}
\end{theorem}

\begin{theorem}\label{T6}
Let $w_1$, $w_2$, $w_3$ be any odd positive integers. Then we have
the following three symmetries in $w_1$, $w_2$, $w_3$:
\begin{equation}\label{N46}
\begin{split}
&(w_1w_2)^n\sum_{a=0}^{w_1d-1}\sum_{b=0}^{w_2d-1}(-1)^{a+b}\chi(ab)
E_{n,\chi}(w_3y_1+\frac{w_3}{w_1}a+\frac{w_3}{w_2}b)\\
=&(w_2w_3)^n\sum_{a=0}^{w_2d-1}\sum_{b=0}^{w_3d-1}(-1)^{a+b}\chi(ab)
E_{n,\chi}(w_1y_1+\frac{w_1}{w_2}a+\frac{w_1}{w_3}b)\\
=&(w_3w_1)^n\sum_{a=0}^{w_3d-1}\sum_{b=0}^{w_1d-1}(-1)^{a+b}\chi(ab)
E_{n,\chi}(w_2y_1+\frac{w_2}{w_3}a+\frac{w_2}{w_1}b).
\end{split}
\end{equation}
\end{theorem}

\begin{theorem}\label{T7}
Let $w_1$, $w_2$, $w_3$ be any positive integers. Then we have the
following two symmetries in $w_1$, $w_2$, $w_3$:
\begin{equation}\label{N47}
\begin{split}
&\sum_{k+l+m=n}\binom{n}{k,l,m}E_{k,\chi}(w_1y)E_{l,\chi}(w_2y)E_{m,\chi}(w_3y)
w_{3}^{k}w_{1}^{l}w_{2}^{m}\\
=&\sum_{k+l+m=n}\binom{n}{k,l,m}E_{k,\chi}(w_1y)E_{l,\chi}(w_3y)E_{m,\chi}(w_2y)
w_{2}^{k}w_{1}^{l}w_{3}^{m}.
\end{split}
\end{equation}
\end{theorem}

\begin{theorem}\label{T8}
Let $w_1$, $w_2$, $w_3$ be any odd positive integers. Then we have
the following two symmetries in $w_1$, $w_2$, $w_3$:
\begin{align}
&\sum_{k+l+m=n}\binom{n}{k,l,m}T_k(w_1d-1,\chi)T_l(w_2d-1,\chi)T_m(w_3d-1,\chi)
w_{3}^{k}w_{1}^{l}w_{2}^{m}\\
&\sum_{k+l+m=n}\binom{n}{k,l,m}T_k(w_1d-1,\chi)T_l(w_3d-1,\chi)T_m(w_2d-1,\chi)
w_{2}^{k}w_{1}^{l}w_{3}^{m}.
\end{align}
\end{theorem}


\end{document}